\documentclass{amsart}
\usepackage{amsfonts}
\usepackage{amsmath}
\usepackage{amssymb}
\usepackage{amscd}
\usepackage{mathrsfs}  
\usepackage{amsbsy}
\usepackage{amsthm}
\usepackage{graphicx}
\usepackage{fancyhdr}
\usepackage{color}

\usepackage[OT2,OT1]{fontenc}  
\newcommand\cyr{%
\renewcommand\rmdefault{wncyr}%
\renewcommand\sfdefault{wncyss}%
\renewcommand\encodingdefault{OT2}%
\normalfont \selectfont} \DeclareTextFontCommand{\textcyr}{\cyr}
\input epsf 

\newcommand{\cD}{\mathcal{D}}

\newcommand{\inn}[2]{\langle\,#1,\,#2\rangle}

\newcommand{\st}{\,|\,}

\newcommand{\R}{{\mathbb{R}}}

\newcommand{\Z}{{\mathbb{Z}}}

\newcommand{\cW}{{\mathcal{W}}}
\newcommand{\wW}{\widetilde{\mathcal{W}}}
\newcommand{\wG}{\widetilde{\Gamma}}
\newcommand{\wK}{\widetilde{K}}

\newcommand{\GLn}{\mathrm{GL}(n,\R)}

\newtheorem{example}{Example}[section]
\newtheorem{remark}{Remark}[section]

\newtheorem{definition}{Definition}[section]

\newtheorem{proposition}{Proposition}[section]
\newtheorem{theorem}{Theorem}[section]
\newtheorem{conjecture}{Conjecture}[section]

\newcommand{\be}{\begin{equation}}
\newcommand{\ee}{\end{equation}}
\newcommand{\nl}{\vskip 5pt\noindent}

\def\sideremark#1{\ifvmode\leavevmode\fi\vadjust{\vbox to0pt{\vss
 \hbox to 0pt{\hskip\hsize\hskip1em
\vbox{\hsize2cm\tiny\raggedright\pretolerance10000
 \noindent #1\hfill}\hss}\vbox to8pt{\vfil}\vss}}}%

\newenvironment{sumlist}
{\begin{list}{$\bullet$}{%
\setlength{\leftmargin}{12.5pt}
\setlength{\labelsep}{6.5pt}
\setlength{\labelwidth}{7.5pt}
\setlength{\parsep}{1pt}
\setlength{\topsep}{0.2\baselineskip}
\setlength{\itemsep}{0.2\baselineskip}}}
{\end{list}}

\begin{document}

\title{Three-way tiling sets in two dimensions}
\author{David Larson, Peter Massopust, and Gestur \'Olafsson}
\address{Department of Mathematics, Texas A{\&}M University, College Station, Texas 77843, USA}
\email{larson@math.tamu.edu}
\address{GSF-National Research Center for Environment and Health,
Institute of Biomathematics and Biometry, and Centre of Mathematics M6,
Technische Universit\"at M{\"u}nchen,
Germany}
\email{massopust@ma.tum.de}
\address{Department of Mathematics, Louisiana State
University, Baton Rouge,
LA 70803, USA} \email{olafsson@math.lsu.edu}
\thanks{The research of the first author was partially supported by NSF
grant DMS-0139386 and of the third author by NSF grant DMS-0402068. The
research of the second author was partially supported by the
grant MEXT-CT-2004-013477, Acronym MAMEBIA, of the European
Commission. The essential part of this research was accomplished
while all three authors participated in the Banff workshop on "Operator
methods in fractal analysis, wavelets and dynamical systems," December 2 -- 7, 2006}
\begin{abstract}

In this article we show that there exist measurable sets $W \subset
\R^2$ with finite measure that tile $\R^2$ in a measurable way under
the action of a expansive matrix $A$, an affine Weyl group $\wW$,
and a full rank lattice $\wG\subset \R^2$.  This note is follow-up
research to the earlier article "Coxeter groups and wavelet sets"
by the first and second authors, and is also relevant to the earlier
article "Coxeter groups, wavelets, multiresolution and sampling" by
M. Dobrescu and the third author. After writing these two articles,
the three authors participated in a workshop at the Banff Center  on
"Operator methods in fractal analysis, wavelets and dynamical
systems," December 2 -- 7, 2006, organized by O. Bratteli, P.
Jorgensen, D. Kribs, G. \'Olafsson, and S. Silvestrov, and discussed
the interrelationships and differences between the articles, and
worked on two open problems posed in the Larson-Massopust article.
We solved part of Problem 2, including a surprising positive
solution to a conjecture that was raised, and we present our results
in this article.
\end{abstract}
\maketitle

\section*{Introduction}
\noindent This article could well have been entitled ``dual wavelet
sets", but we felt that the three-way tiling property was the
important feature to emphasize in the title.  It concerns measurable
sets which are simultaneously dilation-translation and
dilation-reflection wavelet sets, whose existence was shown in
\cite{lm06}, and whose classification has hardly begun. This
represents follow-up research to the article \cite{lm06} and-in
a different way-is related to the article \cite{do05,do06}.

In the article \cite{lm06}, section 7,  the authors show that if
$A\in \GLn$ is an expansive matrix and $\cD=\{A^n\mid n\in\Z\}$ and
$\wW=\cW \ltimes \Gamma$ an affine Weyl group, then there exists a
measurable set $W$ of finite measure such that $\{g W \mid
g\in\cD\}$ and $\{ w W \mid w\in \wW\}$ both form  measurable
tilings of $\R^n$, and moreover such two-way tiling sets always
exist if $A$ is expansive. These were called
$(\cD,\wW)$-\textit{dilation-reflection wavelet sets}. On the other
hand, it is well known that if $\wG\subset \R^n$ is a full rank
lattice, i.e., a co-compact discrete subgroup of $\R^n$, then a
measurable set $W\subset \R^n$ is a
$(\cD,\wG)$-\textit{dilation-translation wavelet set}, in the usual
sense, if and only if $\{A^n W \mid n\in \Z\}$ and $\{W +\gamma\mid
\gamma\in \wG\}$ both form measurable tilings of $\R^n$.  It was
shown in \cite{dls97} that $(\cD,\wG)$-dilation-translation wavelet
sets always exist.    So Theorem 7.4 of \cite{lm06} completely
generalizes Corollary 1 of \cite{dls97} to the case where an
arbitrary affine Weyl group replaces the arbitrary translation group
in the theorem.

In Chapter 8 of \cite{lm06}, the
authors noted that two examples
given of dilation-reflection wavelet sets were
actually "three-way tiling sets" in the sense that they tiled $\R^2$
under three groups: dilation, reflection, and translation. In other
words, they were both dilation-translation wavelet sets in the
traditional sense and also dilation-reflection wavelet sets, for the
same dilation group but with a reflection group (i.e., an affine Weyl
group) replacing the usual translation group. In fact these sets
were known dyadic wavelet sets in the plane: the so-called ``wedding
cake set" and ``four-corners set" given in \cite{dl98}
Not only did they
tile the plane under dilation by $2I$ and translation by $2\pi$
along the $x-y$ coordinate axes, but they also tiled the plane under
an affine Weyl group whose corresponding foldable figure was a
square.

Questions were then raised as to whether three-way-tilers were rare
or common. Our results in the present paper indicate that they are
more common than we expected.  However, the dependence of the full
rank lattice (the translation group) on the foldable figure (and
hence the affine Weyl group) cannot be completely removed. A
simple restriction is that a cell for the lattice must have the same
Lebesgue measure as that of the foldable figure.  The question
is whether this is the only restriction or whether it is also necessary
that the intersection of the lattice with the affine Weyl group
is also a full rank lattice.

We will prove that given any foldable figure $C$ in the plane
containing the origin $0$ in its interior, and given any expansive
dilation matrix $A$, there is a full rank lattice  $\wG\subset \R^n$
which depends on $C$ but not on $A$, such that a measurable set $W$
exists, which tiles $\R^2$ under each of the groups $\wG\subset
\R^n$, $\cD=\{A^n\mid n\in \Z\}$, and $\wW=\cW \ltimes \Gamma$,
where the latter is the affine Weyl group determined by the
reflections about the bounding hyperplanes of $C$.

By applying Proposition 8.1 from \cite{lm06} (see Proposition 3.1
below), the proof reduces to showing that given such a $C$, there is
a full rank lattice $\wW=\cW \ltimes \Gamma$ which has a lattice
cell congruent to $C$ under the action of $\wW=\cW \ltimes \Gamma$
(see definition of congruence below), such that the intersection of
$\wG\subset \R^n$ with the affine Weyl group $\wW=\cW \ltimes
\Gamma$ is itself a full rank lattice.

\section{Preliminaries}
\noindent
In this section we give a short exposition of wavelets and wavelet sets.
We also discuss some relations to finite Coxeter groups and
affine Weyl groups. Standard reference for wavelet sets
are the articles \cite{dls97,dls98,os05,w02} and
for the connection with finite and affine Weyl groups
\cite{lm06,do05,do06}.

\subsection{Translation-dilation wavelet sets}
The $n$-dimensional Fourier transform on $L^2(\mathbb{R}^n)$ is
defined by
\[ (\mathcal{F} f)(s) :=
\int_{\mathbb{R}^n} e^{-2\pi i \langle s , t\rangle} f(t)\,  dm (t) \]
for all $f\in L^1 (\mathbb{R}^n) \cap L^2(\mathbb{R}^n)$. Here $m$ is the
product Lebesgue measure on $\R^n$ and
$\langle \cdot ,\cdot \rangle$ is the
standard inner product. The Fourier transform extends to a
unitary isomorphism of
$L^2(\mathbb{R}^n)$ onto itself.

For $g\in L^1 (\mathbb{R}^n)  \cap L^2(\mathbb{R}^n)$ the
inverse of the Fourier transform is given by
\[ (\mathcal{F}^{-1} g)(t) :=
\int_{\mathbb{R}^n} e^{2\pi i \langle s, t\rangle } g(s) dm \]

Let $A$ be an expansive $n\times n$ real matrix (i.e., all eigenvalues of $A$
have absolute value $> 1$).  (See (\cite{lm06}, Remark 2.1 for six
equivalent definitions of congruence.)  Let $\wG$ be a full rank
lattice in $\mathbb{R}^n$. Equivalently, there exists a basis (not-necessarily orthogonal)
$\{b_1, b_2, ..., b_n\}$  for $\R^n$ such that $\wG$ is the
group of translations by vectors from the group $\{\Sigma k_ib_i \st k_i
\in \mathbb{Z}\}$.  For convenience, let $\wG$ also denote this set
of vectors, i.e., we identify $\wG$ and the subset of $\R^n$
given by $\wG\cdot 0$. By a \emph{dilation - $A$, translation $\wG$ orthonormal
(single) wavelet} we mean a function $\psi \in L^2(\mathbb{R}^n)$
such that
\begin{equation}
\left\{|\det(A)|^{\frac{n}{2}} \psi(A^n t - \gamma) \st n\in \mathbb{Z},
\gamma \in \wG \right\}
\end{equation}
is an orthonormal basis for $L^2(\mathbb{R}^n )$.

Denote by $\chi_E$ the indicator function of a measurable set $E\subseteq \R^n$.
Then $E$ is a \textit{wavelet set} for $A$
and $\wG$ if
\[ \mathcal{F}^{-1} \left(\frac{1}{\sqrt{m(E)}} \,\chi_E\right) \]
is an orthonormal wavelet for $A$ and $\wG$.

\begin{definition} \label{def 7.1}
(\cite{dl98},\cite{dls97}) Let $G$ be a discrete group acting on a
measure space $(M,\mu)$ and let $\Omega, \Sigma\subset M$ be
$\mu$-measurable. Then $\Omega $ and $\Sigma$ are $G$-congruent if
there exists a subset $I\subset G$, a $\mu$-measurable partition
$\{\Omega_g\}_{g\in I}$ of $\Omega$, and a $\mu$-measurable
partition $\{\Sigma_g\}_{g\in I}$ of $\Sigma$ such that
$g\Omega_g=\Sigma_g$.
\end{definition}

It is obvious that if $\Omega$ and $\Sigma$ are
$G$-congruent, then $\Omega$ is a $\mu$- measurable $G$-tile of
$(M,\mu)$  (i.e., a fundamental domain for $G$) if and only if
$\Sigma$ is a $\mu$-measurable $G$-tile for $(M,\mu)$. Moreover, two
different $G$-tiles must be $G$-congruent. If $G$ is a group of
measure-preserving transformations then
$G$-congruence preserves measures of sets. So for such a group, all
$G$-tiles (if any exist) must have the same measure. For a group
which does not consist of measure-preserving transformations, such
as a dilation group on $\R^n$, $G$-congruence can change measures of
sets, and $G$-tiles can have widely differing measures.

A theorem from \cite{dls97} states that a measurable subset $W$ of
$\R^n$ is a dilation-$A$, translation-$\wG$, wavelet set if and only
if $W$ tiles $\R^n$ under both the dilation group $\cD=\{(A^T)^k\mid
k\in \Z\}$  and the translation group
$$\wG^*=\left\{\gamma\in\R^n\mid (\forall \sigma \in\wG)\, \langle \gamma ,\sigma \rangle\in\Z\right\}$$
which is also a full rank lattice, \textit{the dual lattice} of
$\wG$. This was first proven for the case $n = 1$ in \cite{dl98},
and extended to $\R^n$ in \cite{dls97} together with a proof that
such wavelet sets always exist for expansive dilations. It is worth
mentioning here that, under the Fourier transform the translation by
$\gamma\in\wG$ is transformed into modulation by $e^{2\pi i\langle
\gamma ,\cdot \rangle }$. It was first proved by Fuglede in
\cite{f74}, see below, that a measurable set $W$ tiles $\R^n$ under
$\wG^*$ if and only if $\{m(W)^{-1/2}\, e^{2\pi i\langle \sigma
,\cdot \rangle}\mid \sigma \in \wG\}$ forms an orthonormal basis for
$L^2(W)$, which, together with the tiling property under $\cD$ is
needed to construct an orthonormal basis for $L^2(\R^n )$.

\begin{theorem}[Fuglede \cite{f74}]
Assume that $\wG$ is a lattice. Then $L^2(E )$ has an orthogonal basis consisting
of exponentials if
and only if $\{E+t\mid t\in\wG^*\}$ is a measurable tiling of $\mathbb{R}^n$.
\end{theorem}
This result and several examples led Fuglede to conjecture, cf. \cite{f74}:
\begin{conjecture}[The Spectral-Tile Conjecture]
Let $E$ be a measurable subset of $\R^n$. Then $L^2(E)$ has an orthogonal
basis consisting of exponentials  if and only if it is
an additive tile for some discrete subset $\Lambda\subset \R^n$.
\end{conjecture}

Several people worked on this conjecture and derived important results and validated the
conjecture for some special cases, see \cite{ikt99,jp98,lw97,w02}
and the references therein.
However, in 2003,  Tao \cite{t04} showed that the conjecture is false in dimension
$5$ and higher if the lattice hypothesis is
dropped. The other direction was disproved by
Kolountzakis and Matolcsi
in \cite{km06}.
But even now, after  the Spectral-Tiling conjecture has
been proven to fail in higher dimensions, it is still interesting and
important to understand better the connection between spectral properties
and tiling in particular, because of the connection to
wavelet sets.

\subsection{Finite and affine Weyl groups}
In order to proceed, a short excursion into the theory of Coxeter and Weyl groups
as well as foldable figures is necessary. The interested reader is referred to
\cite{Bo,C,G,Gu,H,HW} for more details and proofs.
\subsubsection{Coxeter groups}
\begin{definition}
A Coxeter group $\mathcal{C}$ is a discrete groups with a finite number of
generators $\{r_i\,\st\,i = 1,\ldots, k\}$ satisfying
\[
\mathcal{C} = \bigl\langle r_1,\ldots, r_k\,\st\, (r_i r_j)^{m_{ij}} = 1,\; 1 \leq i,j \leq k\bigr\rangle
\]
where $m_{ii} = 1$, for all $i$, and $m_{ij}\geq 2$, for all $i\neq j$. ($m_{ij}=\infty$ is used to indicate that no relation exists.)
\end{definition}
A geometric representation of a Coxeter group is given by considering it as a subgroup of $GL (V)$, where $V$ is a $k$-dimensional real vector space, which we take to be $\mathbb{R}^k$ endowed with the standard inner product $\inn{\cdot}{\cdot}$. In this representation, the generators are interpreted in the following way.
\par
A reflection about a linear hyperplane $H$ is defined as a linear mapping $\rho: V \to V$ such that $\rho\vert_{H
} = \textrm{id}_H$ and $\rho(x) = -x$, if $x\in H^\perp$. Thus,
$\rho$ is an isometric isomorphism of $V$ of order two such that the multiplicity of the eigenvalue $-1$ is one.
\par
Now suppose that $0\neq r\in H^\perp$. Define
\begin{equation}\label{eq-co-root}
r^\vee :=\frac{2}{\inn{r}{r}}\, r\, .
\end{equation}
Then an easy computation shows that
\begin{equation}\label{refl}
\rho_r(x) = x - \inn{x}{r^\vee}\, r =x -\inn{x}{r}\, r^\vee
\end{equation}
is the reflection about the hyperplane $H$ perpendicular to $r$.
\par
If $R\subset \R^n$ is a finite set such that the group $<\rho_r\mid r\in R>$ generated by the reflection
$\rho_r$, $r\in R$, is finite, then it is a finite Coxeter group. In particular this is the
case if $R$ is a root system.
%
\subsubsection{Roots systems and Weyl groups}
The normal vectors to a set of hyperplanes play an important role in the representation theory for Coxeter groups. We have seen above that they
correspond to the generators of such groups. Two such normal vectors, $\pm\, r$, that are orthogonal to a hyperplane are called {\em roots}.
\begin{definition}
A \textit{root system} ${R}$ is a finite set of nonzero vectors $r_1, \ldots, r_k\in$ $\mathbb{R}^n$ satisfying
\begin{enumerate}
\item   $\mathbb{R}^n = \textrm{span}\,\{r_1,\ldots, r_k\}$
\item   $r, \alpha r\in R$ if and only if $\alpha = \pm 1$
\item   $\forall r\in R$: $\rho_r(R)= R$, i.e., the root system $R$ is
        closed with respect to the reflection through the hyperplane orthogonal to $r$.
\item   $\forall r,s\in R$: $\displaystyle{\inn{s}{r^\vee}\in\mathbb{Z}}$.
\end{enumerate}
Note that (4) is a strong restriction on $\theta=\measuredangle (r,s)$, the angle between
$r$ and $s$. It implies that
$$4\cos (\theta )^2=4\frac{\inn{r}{s}^2}{\|r\|^2 \| s\|^2}\in \Z\, .$$
\end{definition}

Assume that $R$ is a root system. Then the group generated by
$\rho_r$, $r\in R$, is a finite  Coxeter group. It is called the \textit{Weyl Group} $\mathcal{W}$ of $R$.

A subset $R^+\subset R$ is called a set of \textit{positive} roots if
$R^++R^+\subseteq R^+$, $R^+\cap (-R^+)=\emptyset $ and
$R=R^+\cup -R^+$. Let $v\in\mathbb{R}^n$ such that $\inn{r}{v}\not= 0$ for all $r\in R$.
Then the set $R^+=\{r\in R\mid \inn{r}{v}>0\}$ is a positive system and all positive
system can be constructed in this way. Roots that are not positive are called
\textit{negative}.
\par
\begin{example}
A simple example of a Weyl group in $\mathbb{R}^2$ is given by the root system depicted in Figure \ref{klein}. The roots are $r_1 = -r_3 =
(1,0)^\top$ and $r_2 = - r_4 = (0,1)^\top$. The positive roots are $r_1$ and $r_2$. The group of reflections generated by these four roots
is given by
\[
V_4 := \bigl\langle \rho_1,\rho_2\,\st\, \rho_1^2 = \rho_2^2 = 1,\; (\rho_1\rho_2)^2 =1\bigr\rangle,
\]
where $\rho_1$ and $\rho_2$ denotes the reflection about the $y$-, respectively, $x$-axis. This group is commutative and called {\em Klein's four-group} or the {\em  group of order four}. In the classification scheme of Weyl groups $V_4$ is referred to as $A_1\times A_1$ since it is the direct product of the group $A_1 := \bigl\langle \rho_1\,\st\, \rho_1^2 = 1\bigr\rangle$ whose root system is $R = \{r_1, r_3\}$ with itself.
\begin{figure}[h]
\begin{center}
\includegraphics[width=1in,height=1in]{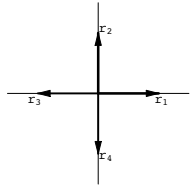}
\end{center}
\caption{The root system for Klein's four-group.}\label{klein}
\end{figure}
\end{example}
\par
For the following, we need some properties of roots systems and Weyl groups, which we state in a theorem.
\begin{theorem}
Let $R$ be a root system and $\mathcal{W}$ the associated Weyl group. Then the following hold.
\begin{enumerate}
\item   Every root system $R$ has a basis $\mathcal{B} = \{b_i\}$ consisting of positive (negative) roots.
\item   Let $C_i := \{x\in\mathbb{R}^n\,\st\,\inn{x}{b_i} > 0\}$ be the {\em Weyl chamber} corresponding to the basis $\mathcal{B}$. Then
        the Weyl group $\mathcal{W}$ acts simply transitively on the Weyl chambers.
\item   The set $C := \overline{\bigcap_i C_i}$ is a noncompact fundamental domain for the Weyl group $\mathcal{W}$. It is a simplicial cone, hence
        convex and connected.
\end{enumerate}
\end{theorem}
In order to introduce foldable figures below, we need to consider reflections about affine hyperplanes. For this purpose, let $R$ be a root
system. An {\em affine hyperplane} with respect to $R$ is given by
\begin{equation}\label{hyper}
H_{r,k} := \{x\in\mathbb{R}^n\st\inn{x}{r} = k\},\qquad k\in\mathbb{Z}.
\end{equation}
It is easy to show that reflections about affine hyperplanes have the form
\begin{equation}\label{affref}
\rho_{r,k}(x) = x - \displaystyle{\frac{2(\inn{x}{r}-k)}{\inn{r}{r}}}\,r =: \rho_r (x) + k\,{r}^\vee,
\end{equation}
where $r^\vee := 2\,r/\inn{r}{r}$ is the {\em coroot} of $r$.
\begin{definition}
The {\em affine Weyl group} $\widetilde{\mathcal{W}}$ for a root system $R$ is the (infinite) group generated by the reflections $\rho_{r,k}$ about the affine hyperplanes $H_{r,k}$:
\[
\widetilde{\mathcal{W}} := \bigl\langle \rho_{r,k}\st r\in R, k\in\mathbb{Z}\bigr\rangle
\]
We sometimes  will refer to the concatenation of elements from $\widetilde{\mathcal{W}}$ as {\em words}.
\end{definition}
\begin{theorem}\label{th4.5}
The affine Weyl group $\widetilde{\mathcal{W}}$ of a root system $R$ is the semi-direct product $\mathcal{W}\ltimes \Gamma$, where $\Gamma$
is the abelian group generated by the coroots ${r}^\vee$. Moreover, $\Gamma$ is the subgroup of translations of $\widetilde{\mathcal{W}}$ and
$\mathcal{W}$ the isotropy group (stabilizer) of the origin. The group $\mathcal{W}$ is finite and $\Gamma$ infinite.
\end{theorem}
\begin{remark}
There exists a complete classification of all irreducible affine Weyl groups and their associated fundamental domains. These groups are given as
types $A_n$ ($n\geq 1$), $B_n$ ($n\geq 2$), $C_n$ ($n\geq 3$), and $D_n$, ($n\geq 4$), as well as $E_n$, $n = 6, 7, 8$, $F_4$, and $G_2$. (For more
details, we refer the reader to \cite{Bo} or \cite{H}.)
\end{remark}
We need a few more definitions and related results. By a \textit{reflection group} we mean a group of transformations generated by the reflections
about a finite family of affine hyperplanes. Coxeter groups and affine Weyl groups are examples of reflections groups.
\par
Let $\mathcal{G}$ be a reflection group and $\mathcal{O}_n$ the group of linear isometries of $\mathbb{R}^n$. Then there exists a homomorphism
$\phi: \mathcal{G}\to\mathcal{O}_n$ given by
\[
\phi(g)(x) = g(x) - g(0),\quad g\in\mathcal{G},\; x\in\mathbb{R}^n.
\]
The group $\mathcal{G}$ is called essential if $\phi(\mathcal{G})$ only fixes $0\in\mathbb{R}^n$. The elements of $\ker\phi$ are called translations.
\subsection{Foldable figures}
\begin{figure}[h]
\begin{center}
\includegraphics[width=2cm,height=2cm]{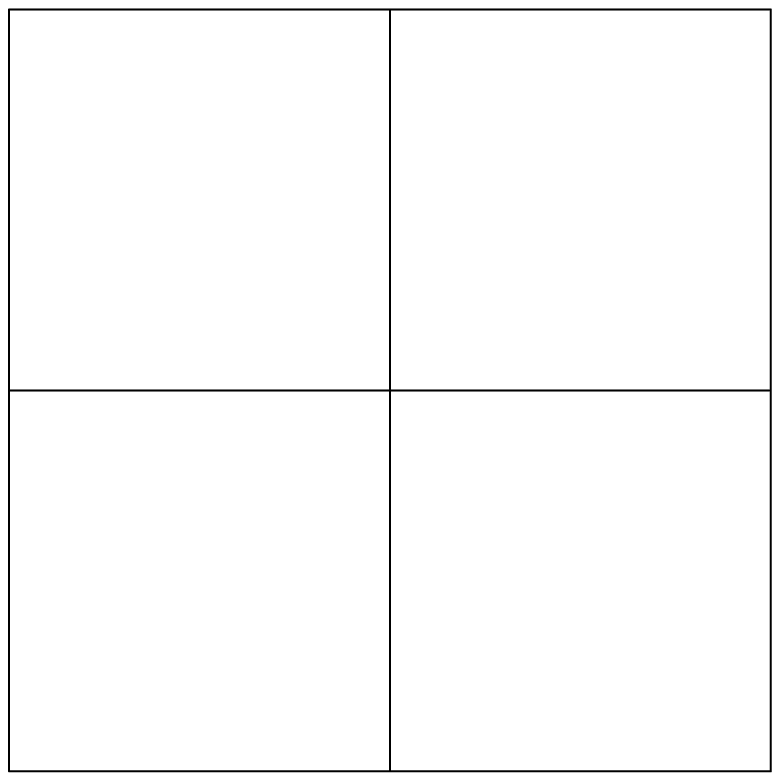}\hspace{2cm}
\includegraphics[width=2cm,height=2cm]{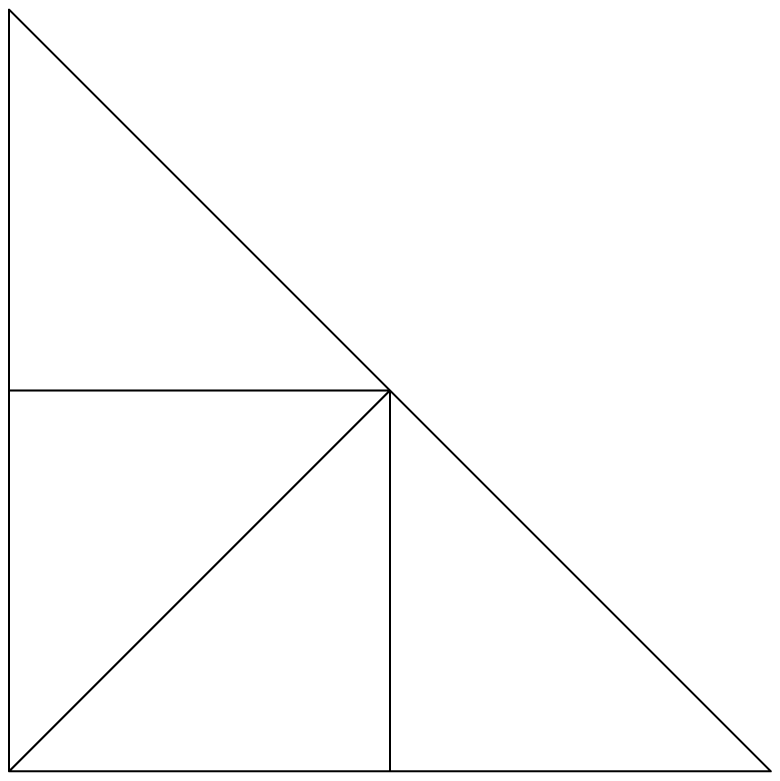}
\end{center}
\caption{Examples of foldable figures.}\label{fig5}
\end{figure}

In this subsection, we define for our later purposes the important concept of a foldable figure \cite{HW}.
\begin{definition}
A compact connected subset $F$ of $\mathbb{R}^n$ is called a \textit{foldable figure} if and only if there exists a finite set $\mathcal{S}$ of affine
hyperplanes that cuts $F$ into finitely many congruent subfigures $F_1, \ldots, F_m$, each similar to $F$, so that reflection in any of the cutting
hyperplanes in $\mathcal{S}$ bounding $F_k$ takes it into some $F_\ell$.
\end{definition}
In Figure~\ref{fig5} are two examples of foldable figures shown.
Properties of foldable figures are summarized in the theorem below. The statements and their proofs can be found in \cite{Bo} and \cite{HW}.
\begin{theorem}
${}$
\begin{enumerate}
\item   The reflection group generated by the reflections about the bounding hyperplanes of a foldable figure $F$ is the affine Weyl group
        $\widetilde{W}$ of some root system. Moreover, $\widetilde{W}$ has $F$ as a fundamental domain.
\item   Let $\mathcal{G}$ be a reflection group that is essential and without fixed points. Then $\mathcal{G}$ has a compact fundamental domain.
\item   There exists a one-to-one correspondence between foldable figures and reflection groups that are essential and without fixed points.
\end{enumerate}
\end{theorem}

\subsection{Tiling sets for dilations and affine Weyl groups}
Recall, that a foldable figure is a connected, convex fundamental
domain for an affine Weyl group $\wW$.

\begin{definition} \label{def7.2} (Def 7.3 from \cite{lm06})
Given an affine Weyl group $\widetilde{{\mathcal{W}}}$ acting on
$\mathbb{R}^n$ with fundamental domain a foldable figure $C$, given
a designated interior point $\theta$ of $C$, and given an expansive
matrix $A$ on $\mathbb{R}^n$, a dilation-reflection wavelet set for
$(\widetilde{{\mathcal{W}}}, \theta, A)$ is a measurable subset $W$
of $\mathbb{R}^n$ satisfying the properties:
\begin{enumerate}
\item $W$ is congruent to $C$ (in the sense of Definition 2.4) under the
action of $\widetilde{{\mathcal{W}}}$, and
\item $W$ generates a measurable partition of $\mathbb{R}^n$ under the action of
the affine mapping $D(x) := A(x - \theta) + \theta$.
\end{enumerate}
In the case  where $\theta = 0$, we abbreviate
$(\widetilde{{\mathcal{W}}}, \theta, A)$ to
$(\widetilde{{\mathcal{W}}}, A)$.
\end{definition}

In the above definition, the reason $\theta$ is prescribed is that
in order to apply the methods of \cite{dls97} for two-way tilers, it
is simplest to assume that the dilation fixed point be contained in
the interior of the model for abstract translation tiler, namely the
foldable figure $C$ in this case. (It is not a necessary
restriction, but on the other hand some kind of restriction such as
this is necessary. It is a sufficient restriction to make the
apparatus work.) Then, to make the dilation by $A$ \emph{compatable}
with the action of the reflection group (i.e., so that $\theta$ is a
fixed point for the dilation transformation group) the usual
dilation by $A$ is replaced with the affine dilation mapping  $D(x)
:= A(x - \theta) + \theta$. In practice, in proofs in particular, we
frequently perform an initial translation of the system so that
$\theta = 0$, and dilations thus correspond to the usual dilation
group. On the other hand, in reflection group theory one usually
performs an initial translation of a geometric problem so that the
orgin $0$ lies at a \emph{vertex} of the foldable figure, hence not
an interior point. Hence, in particular in the analysis of the
\emph{root systems} (see below), this is the assumption made.  A
simple translation of the system then yields the general case. Thus,
basically, in order to apply our theory (in its present form), we
either need to assume up-front that $0$ is an interior point of $C$
and use the usual dilation group $\cD = \{A^k W \mid n\in \Z\}$, or
we need to fix a designated point $\theta$ in the interior of $C$
and use the affine dilation group generated by $D(x) := A(x -
\theta) + \theta$.

As a foldable figure $C$ supports an orthonormal basis for
$L^2(C)$ of fractal surface functions, similar to the way a
dilation-translation wavelet set $W$ supports an orthonormal basis
of exponential functions, and if we compose these fractal surface
functions with the reflection group congruence operations we obtain
an orthonormal basis of fractal-surface-induced functions on the
dilation-reflection wavelet set $W$.  Since $W$ tiles $\R^n$ under
dilation by $A$, by the dilations by powers of $A$ to these
functions we obtain a fractal surface induced orthonormal basis of
$L^2(\R^n)$.  So in a function-theoretic sense as well as in a set
theoretic sense (i.e., tiling properties), the dilation-reflection
wavelet sets are natural generalizations of the dilation-translation
wavelet sets: they are both affiliated with orthonormal bases of
$L^2(\R^n)$.

The main result of \cite{lm06} is:

\begin{theorem}[Theorem 7.4 from \cite{lm06}]\label{100}
There exist $(\widetilde{{\mathcal{W}}}, \theta, A)$-wavelet sets
for every choice of $\widetilde{W}$, $\theta$, and $A$.
\end{theorem}

\subsection{Subspace wavelet sets and Weyl groups}
In \cite{do05,do06} the connection between
Weyl groups (of finite reflection groups) was
looked at from a slightly different point of view.

Most of the examples of wavelet sets tend to be fractal like and symmetric around $0$.
The first aim was to construct wavelet functions that
have some directional properties in
the frequency domain. The direction is given by a fundamental domain $C$ for
the Weyl group, which is also a convex cone, i.e., $C+C\subset C$ and
$\R^+\, C\subset C$. The wavelets are again given
by the inverse Fourier transform of the normalized indicator function
of a measurable subset $E$. Assuming that $A^T(C)\subset C$, and $A^T\wG^*\subset \wG^*$, then it was shown, that
one can find a set $E\subset C$, such that $E$ tiles $C$ under $\{(A^T)^k\mid k\in\Z\}$
and tiles $\R^n$ under $\wG^*$, i.e, $E$ is a
$(A,\wG)$-subspace wavelet set. The Weyl group is then
used to rotate the frequency domain to cover all of $\R^n$.
The construction was still fractal, but the result  is not symmetric around $0$ anymore.

A natural question one asks  when working with wavelets is whether they are related
with any multiresolution analysis or multiwavelets.
It was shown that it is possible to construct subspace
multiwavelets using the Weyl group. In particular, the wavelets are still directional in
the frequency domain. In fact, the support
of the Fourier transform is supported in cones which are
fundamental domains for the action of a Weyl group on the Euclidian space.
Finally, the  relation between those constructions and  sampling theory
were discussed.
Any square-integrable function can be written as a sum of its projections on subspaces,
where each subspace contains only signals supported in the frequency
domain in the cones mentioned above.
Each projection can then be sampled using a version of the
Whittaker-Shannon-Kotel'nikov sampling theorem.

\section{Three-way tiling sets}
\noindent
The following proposition was proven in \cite{lm06}, which covers
the cases of the four-corners and wedding-cake sets, as mentioned in
the introduction.  We will give a less general version than we gave
in \cite{lm06} because it suffices for our purposes in this article.

\begin{proposition}
[Proposition 8.1 in \cite{lm06}]
 Suppose that $C$ is any
foldable figure which is a fundamental domain for both a translation
group $\wG$ and the affine Weyl group $\wW$ for $C$, and which
contains $0$ in its interior. If the intersection group
$\mathcal{J}$ of $\wG$ and $\wW$ contains a full rank lattice, then
for any expansive matrix $A$ there exist sets $W$ which are
simultaneously dilation $A$-translation $\wG$ and dilation $A$-reflection $\wW$ wavelet sets.
\end{proposition}

The version of the above proposition that was stated in \cite{lm06}
only required that the intersection group $\mathcal{J}$ in the
hypothesis be large enough so that the prescribed dilation group
$\cD$ together with $\mathcal{J}$ form an abstract
dilation-translation pair in the sense of \cite{dls97}. (See also
\cite{lm06}, Definition 2.8). This hypothesis is automatically
satisfied if  $\mathcal{J}$ contains a full rank lattice, which is
what we show in the present paper.

In \cite{lm06}, it was not clear whether the above proposition could
be applied to decide whether three-way-tilers could exist in greater
generality than the special cases worked out in that paper. We quote
a problem that the authors of \cite{lm06} pose at the end of section
8. The main purpose of this paper is to provide some examples that
settle a part of this problem affirmatively, indicating that
three-way-tilers could be common. This is the part we state below as
Problem 1a.
\medskip

\noindent PROBLEM 1 (This is Problem 2 from \cite{lm06}): Let $C$ be
any foldable figure in $\R^n$ containing 0 in its interior and let
$\wW=\mathcal{W} \ltimes \Gamma$ be the associated affine Weyl
group. Suppose that $A$ is any expansive matrix in $\mathrm{M}_n (\R
)$, and $\wG$ a full rank lattice in $\R^n$. Let $[0,b_1)\times
[0,b_2)\times \ldots \times [0,b_n)$ be a fundamental domain for
$\wG$. Give necessary and sufficient conditions for the existence of
a set $W$ which is simultaneously

\begin{enumerate}
\item $\wW$-congruent to $C$;
\item a $\cD=\{A^n\mid n\in \Z\}$ measurable tiling set of $\R^n$;
\item $\wG$-congruent to the set $[0,b_1) \times [0,b_2) \times
\ldots \times [0,b_n)$, i.e., a $\wG$-spectral set.
\end{enumerate}

It was noted in \cite{lm06} that any $W$ satisfying (1), (2), and
(3) would be both a dilation-translation wavelet set for $(\cD, \wG
)$ and a dilation-reflection wavelet set for $(\wW,A)$, and
conversely, any set which is both a dilation-translation wavelet set
for $(\cD, \wG )$ and a dilation-reflection wavelet set for $(W,A)$
must satisfy (1), (2), and (3). In particular, the question was
asked:
\medskip

\noindent
PROBLEM 1a. Does there exist such a $W$ for an irreducible affine
Weyl group, such as the group corresponding to a equilateral
triangle, which is a foldable figure? We wrote: \textit{``We think
that the answer is probably no. But in the topic of wavelet sets
there are often surprises, so we would not be very surprised if the
answer was yes.''}

As explained in the introduction, a main point to this paper is to
show that indeed such sets $W$ exist for the equiangular triangle
example and other foldable figures.  We give concrete examples, and
pose some further questions.

We analyze Problem 1 above. This problem has two distinct parts.  To
describe them more clearly, let us define a triple $(\wW, \wG, \cD)$
as above to be \textit{allowable} if there exists a measurable set $W$
satisfying (1), (2), (3)  of Problem 1.

The first part of Problem 1 asks for necessary and sufficient
conditions for such a triple to be allowable.  If a triple is
allowable it follows that the foldable figure for $\wW$ must have
the same volume (measure) as a lattice cell (parallelopiped) for
$\wG$, so not every triple is allowable. In this case we must have
$\wG =B\Z^n$ with $|\det B |=|C|$. But we know of no other general
restriction that we can prove.

Attempting to address this situation, the second part of the
problem, which we list as Problem 1a, asks if the affine Weyl group
in an allowable triple can ever be irreducible, and asks
specifically if the affine Weyl group for Fig. 10 (in \cite{lm06}),
i.e., an equilateral triangle, can be a member of an allowable
triple.  This is the simplest irreducible $\wW$. The authors in
\cite{lm06} conjectured ``no", because all the examples they worked
out had reducible Weyl groups and the simplest irreducible case
seemed intractable.

But in Banff, the present authors showed that the answer to the Fig. 10
problem is actually ``yes", showing that there are still in fact
surprises out there in the theory of wavelet sets, and this surprise
led to the rest of our work in this article. We show that in two
dimensions, given any foldable figure with associated affine Weyl
group $\wW$, and given any expansive matrix $A$, there always exists
a full rank lattice $\wG$ so that the triple $(\wW, \wG, \cD)$ is
allowable.

The key is to construct $\wG$ so that the intersection of $\wG$ with
$\wW$ contains the translation group of some full rank lattice, and
then we apply Proposition 8.1 in \cite{lm06}.  This idea was
outlined and used in \cite{lm06}, but except for certain reducible
cases we did not know how to construct $\wG$ from the foldable
figure. This construction is what was worked out in Banff, and is
the essence of this paper.

We use $\wW$ for an affine Weyl group: $\wW =
\mathcal{W}\ltimes \Gamma$, where $\mathcal{W}$ is the stabilizer of
the origin and $\Gamma$ the translation group generated by the
coroots. Note that $\cW$ is a finite Coxeter group. In the
following, we are abusing language by identifying the (discrete)
abelian group $\Gamma$ with the (geometric) lattice that is
generated by it, i.e., $\Gamma \simeq \wW\cdot 0$. The system of
roots corresponding to $\cW$ is denoted by $R$ and the system of
coroots by $R\,\check{}$. There are only four rank 2 root systems:
$A_1\times A_1$ (reducible with fundamental domain
$[0,1]\times[0,1]$), $A_2$, $B_2$, and $G_2$.
 We will only discuss the root system $A_2$ in
all details, and only list the necessary information and ideas for
the other two irreducible systems.  Our notation will be the same as
in \cite{lm06}.

For the construction of a dilation tiling set one needs that the
fixed point $0$ is in the interior of the foldable figure $C$. But
to simplify the discussion in the following, we assume that $0$ is
one of the vertices of $C$.

Let $e_1 := (1,0)^\top$ and $e_2 := (0,1)^\top$ denote the vectors
of the standard $\R^2$-basis and -- again by abuse of notation --
we will also write $e_1 := (1,0,0)^\top$, $e_2 := (0,1,0)^\top$, and
$e_3 := (0,0,1)^\top$ for the vectors of the standard basis of
$\R^3$. We denote by $L$ the $\Z$-module/lattice $L := \Z
e_1\oplus\Z e_2\oplus\Z e_3$. The cross product in $\R^n$ is denoted
by $\wedge$ and the Euclidean length by $\|\,\cdot\,\|$.

\begin{theorem} Let $R$ be an irreducible root system in $\mathbb{R}^2$ and let $\wW =
\mathcal{W}\ltimes \Gamma$ be the associated affine Weyl group. Then
there exists a full rank lattice $\wG$ and a measurable set
$\Omega\subset \mathbb{R}^2$ such that

\begin{enumerate}

\item $\wG\cap \wW= \wG\cap \Gamma$ is a full rank lattice in $\mathbb{R}^2$;

\item $\Omega$ is a tiling set for $\wG$ and $\wW$.

 \end{enumerate}
 \end{theorem}
\begin{proof} The proof will be done by case by case inspection of
the three root systems $A_2, B_2$, and $G_2$.  Since every case in
$R^2$ is equivalent to one of these, the proof will be complete.  In
each case we will explicitly construct $\wG$ and $\Omega\subset
\mathbb{R}^2$,  and show that they have the properties required. For
a given root system, the lattice $\wG$ and the set $\Omega\subset
\mathbb{R}^2$ are not unique, and no general construction of
examples is known.  (This is an open direction.)

\subsection{Root System $A_2$}\label{s-1}
\noindent In this section we discuss the root system $A_2$. Let
$V\subset\R^3$ be the hyperplane given by the linear equation $x + y
+ z = 0$. We identify $V$ with $\R^2$ by using $\alpha =e_1-e_2$ and
$\beta =e_2-e_3$ as basis for the two dimensional subspace $V$.

The root system $A_2$ is given by:
\begin{eqnarray*}
A_2&=&\{\alpha\in L\cap V\mid \|\alpha \|^2=2\}\\
&=&\pm \{e_1-e_2,e_2-e_3,e_1-e_3\}\, .
\end{eqnarray*}
It is clear from this that $A_2^\vee = A_2$. Note that
$\{e_1-e_2,e_2-e_3,e_1-e_3\}$ is a system of positive roots, with
$\{\alpha,\beta\}$ as a corresponding set of simple roots.
\par
Since $\alpha$ and $\beta$ are in $V$, the angle $\theta$ between
them can be computed via $-1 = \inn{\alpha}{\beta}  = \|\alpha\|\,
\|\beta\|\,\cos \theta$. This yields $\theta = 2\pi/3$.
\par
Now introduce a coordinate system $(v_1,v_2)$ in $V$, with the
$v_1$-axis along the vector $\alpha$ and the $v_2$-axis
perpendicular to it. Since $\alpha$ has length $\sqrt{2}$ its
coordinates with respect to the $(v_1,v_2)$-coordinate system are
$\alpha = (\sqrt{2},0)^\top$. Then the representation of $\beta$ in
the $(v_1,v_2)$-coordinate system is obtained by rotating $\alpha$
by $\theta$ counterclockwise. This gives for $\beta$:
\[
\beta = \begin{pmatrix} \cos 2\pi/3 & -\sin 2\pi/3 \\ \sin 2\pi/3 &
\cos 2\pi/3\end{pmatrix}\begin{pmatrix} \sqrt{2}\\ 0\end{pmatrix} =
\begin{pmatrix} -\sqrt{2}/2\\ \sqrt{6}/2\end{pmatrix}
\]
The third positive root is $\alpha+\beta = (\sqrt{2}/2,
\sqrt{6}/2)^\top$.
\par
Therefore, since the roots are identical to the coroots, the full
rank lattice $\Gamma$ is generated by the simple roots $\alpha$ and
$\beta$:
\[
\Gamma = \Z\alpha \oplus\Z\beta.
\]
The walls of the fundamental domain $C$ for $\wW$ are given by the
linear hyperplanes corresponding to the simple roots, i.e.,
$H_\alpha$ and $H_\beta$ and the affine hyperplane
$H_{\tilde{\alpha},1}$ corresponding to the highest root
$\tilde{\alpha} = \alpha + \beta$. A quick calculation shows that
$H_{\tilde{\alpha},1}$ is represented by the linear equation (in the
$(v_1,v_2)$-coordinate system)
\[
\frac{1}{\sqrt{3}}v_1 + v_2 = \frac{\sqrt{6}}{3}.
\]
Thus, $C$ is a $(\pi/3,\pi/3,\pi/3)$-triangle with vertices at
$(0,0)$, $(\sqrt{2}/2,\sqrt{6}/6)$, and $(0,\sqrt{6}/3)$. \nl We
recall that the affine Weyl group $\wW$ is generated by the
reflections about the walls of $C$. The area of a fundamental
lattice cell $K$ of $\Gamma$ and the area of $C$ are easily computed
to be:
\[
|K| = \sqrt{3}\qquad\textrm{and}\qquad |C| = \sqrt{3}/6.
\]
Note that $|K| = |\cW|\cdot |C|$. The reason for this equality lies
in the fact that for an affine Weyl group $\wW = \cW \ltimes \Gamma$
with fundamental domain $C$, the set $\cW \cdot C$ is a fundamental
domain for the translation lattice. As the Weyl group acts freely on
the positive Weyl chamber, and $C$ can be taken inside a fixed Wey
chamber, it follows that  $|\cW\cdot C| = |\cW| |C|$.
\\
\begin{figure}[h]
\begin{center}
\includegraphics[width=2in,height=2in]{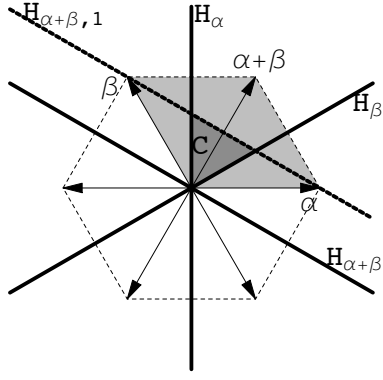}
\caption{Geometry for $A_2$.}\label{f1a}
\end{center}
\end{figure}
\par
In order to obtain three-way tiling sets, we need to find a set
$\Omega$ with the property that $\Omega$ transforms to $C$ via the
action of $\wW$ and  to a fundamental cell $\wK$ via the action of
$\wG$.

For the root system $A_2$, we will exhibit two such sets, $\Omega_1$
and $\Omega_2$. To obtain the first, define
\[
\delta_1 := \frac{1}{6}[(\alpha+\beta) + \beta)] = (0, \sqrt{6}/6)^\top,\quad \eta_1 := \frac{\alpha}{2},
\]
and
\[
\wG_1 := \Z\eta_1\oplus\Z\delta_1.
\]
Observe that $6\,\wG_1\subset\Gamma$ and hence $\wG_1\cap \wW=\wG_1\cap \Gamma$ is a full rank lattice. Furthermore $\delta_1\in\wG_1\setminus\Gamma$ so $\wG_1\not=\Gamma$.

Denote by $\wK_1 := \eta_1\wedge\delta_1$ a fundamental cell
in $\wG_1$, and set
\[
\wK_{11} := \left\{(u,v)\in\wK_1\st v \geq u/\sqrt{3}\right\}\qquad\textrm{and}\qquad\wK_{12} := \wK_1
\setminus\wK_{11}.
\]
Now take as $\Omega_1$ the set
\[
\Omega_1 := \wK_{11}\cup (\wK_{12} + \eta_1 + \delta_1).
\]
By construction, $\Omega_1$ is $\wG_1$-congruent to the fundamental cell $\wK_1$.
Moreover, $\Omega_1$ is also $\wW$-congruent to $C$, for
\[
C = \wK_{11} \cup\left[ \rho_\alpha (\wK_{12} + \eta_1 + \delta_1) + 2 \eta_1\right],
\]
where $\rho_\alpha\in\wW$ denotes the reflection about the linear hyperplane $H_\alpha$. (Note that
$\alpha\in\Gamma$.)
\par
The second set $\Omega_2$ is obtained by taking, as above, $\delta_2 = \frac{1}{6}(\alpha + 2 \beta) = (0, \sqrt{6}/6)^\top$, and defining
\[
\wG_2 := \Z\, \delta_2\oplus\Z\, \eta_2,
\]
where $\eta_2 := \frac{1}{3}(2 \alpha + \beta) = (\sqrt{2}/2, \sqrt{6}/6)^\top$. Note
that $\wG_2\cap\Gamma\supset 6 \wG_2$ and hence $\wG_2\cap \Gamma =
\wG_2\cap \wW$ is a full rank lattice. The fundamental cell $\wK_2
:= \delta_2\wedge\eta_2$ of $\wG_2$ is partitioned into
\[
\wK_{21} := \left\{(u,v)\in\wK_2 \st v \leq -\frac{1}{\sqrt{3}} u + \frac{\sqrt{6}}{3}\right\}\qquad\textrm{and}\qquad\wK_{22} := \wK_2 \setminus\wK_{21}.
\]
Now define \[ \Omega_2 := \wK_{21}\cup (\wK_{22} + \eta_2 - \delta_2). \] Then $\Omega_2$ is $\wG$-congruent to
$\wK_2$ and also $\wW$-congruent to $C$, since
\[
C = \wK_{21}\cup \left[\rho_\alpha (\wK_{22} + \eta_2 - \delta_2) + \alpha\right],
\]
with $\rho_\alpha\in\wW$ denoting the reflection about the hyperplane $H_\alpha$.
\begin{figure}[h]
\begin{center}
\includegraphics[width=4.5cm,height=4.5cm]{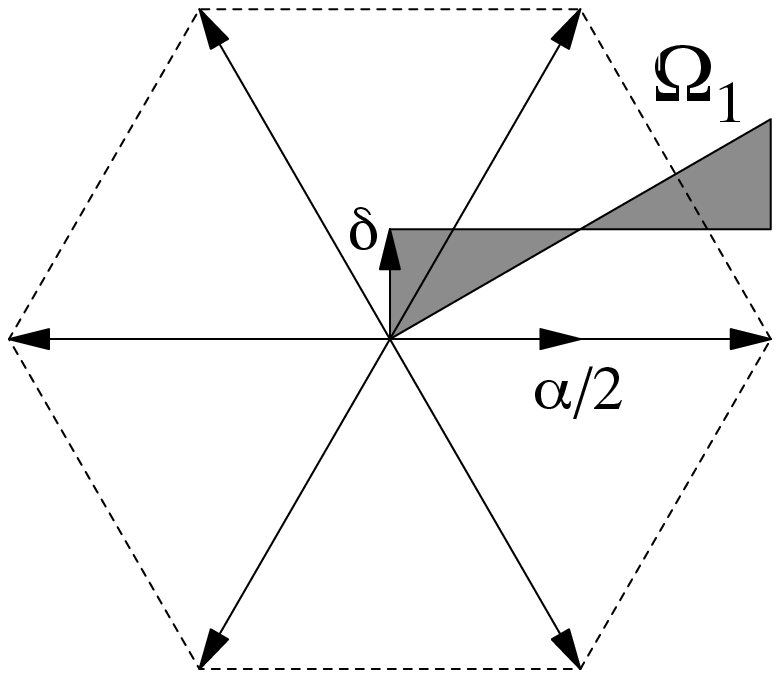}\hspace*{1cm}
\includegraphics[width=4.5cm,height=4.5cm]{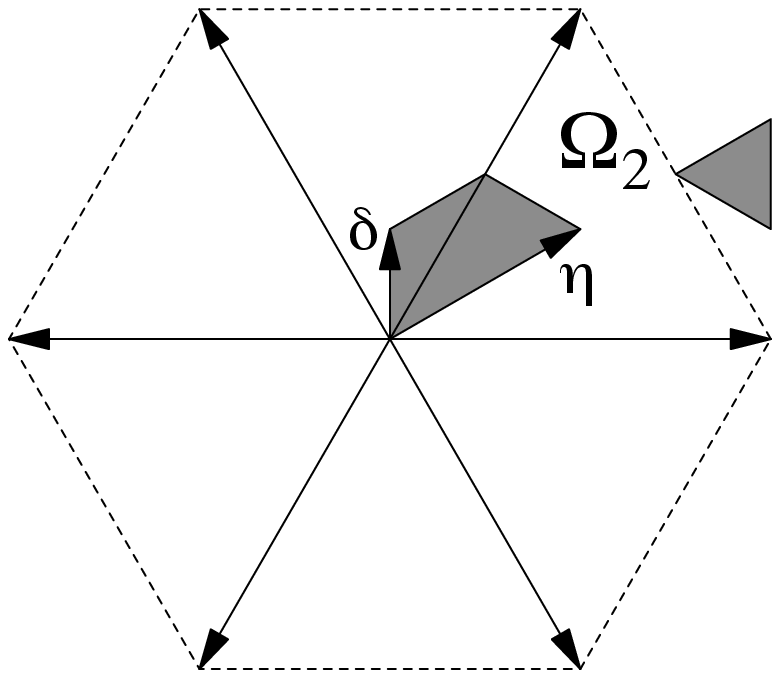}
\caption{The sets $\Omega_1$ and $\Omega_2$ for $A_2$.}\label{f12}
\end{center}
\end{figure}
\subsection{Root System $B_2$}
\noindent In this section we discuss the root system $B_2$. The
arguments are quite similar to the ones for $A_2$ so we restrict our
exposition to list the objects that are needed and then give a some
more detailed arguments at the end.
\\
\begin{sumlist}
\item Root system $R$: $\pm e_1$, $\pm e_2$, $\pm e_1 \pm e_2$.

\item Positive roots: $e_1$, $e_2$, $e_1\pm e_2$.
\item Basis for $\R^2$: $\alpha := e_1 - e_2$ and $\beta := e_2$.
\item Coroot system $R^\vee$:  $\pm 2 e_1$, $\pm 2 e_2$, $\pm e_1 \pm e_2$.

\item Highest root: $\tilde{\alpha} = e_1 + e_2$.
\item Weyl group $\cW$:  $\frak{S}_2$ acting as permutations of the coordinates, and all the sign changes
$x_j\mapsto \pm x_j$. The Weyl group has therefore order $|\cW |=
8$.
\item  Lattice  generated by the coroots: $\Gamma = \Z\alpha{}^\vee\oplus \Z \tilde{\alpha}$.
\item Walls of fundamental domain: $H_\alpha: x - y = 0$, $H_\beta: y = 0$, and $H_{\tilde{\alpha},1}: x + y = 1$.

\item Fundamental domain for $\widetilde{\mathcal{W}}$:
$$C = \{(x,y\in\R^2 \st 0 \leq x+y \leq 1,\; 0\leq y \leq x\}\, .$$
$C$ is a triangle with vertices $(0,0)$, $(1/2,1/2)$, and $(1,0)$,
hence has angles $\pi/4$, $\pi/2$, $\pi/4$.
\item $|K| = 2$ and $|C| = 1/4$.
\end{sumlist}
\par
\begin{figure}[h]
\begin{center}
\includegraphics[width=2in,height=2in]{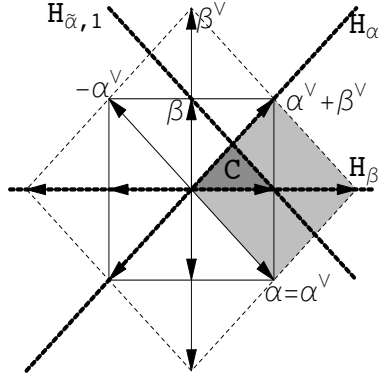} \caption{Geometry for $B_2$.}\label{f2}
\end{center}
\end{figure}
\par
Now choose vectors $\delta_1$ and $\eta_1$ as follows: $ \delta_1 := \frac{1}{2}\tilde{\alpha} = (\frac{1}{2},\frac{1}{2})^\top$ and $\eta_1 := \frac{1}{2}(\alpha+\beta) = (\frac{1}{2},0)^\top$. Let
\[
\wG_1 := \Z\delta_1\oplus\Z\eta_1.
\]
Then, $\wG_1\cap\wW=\wG_1\cap\Gamma \supseteq 4\wG_1$, which is a full rank lattice. Hence
$\wG_1\cap\wW=\wG_1\cap\Gamma $ is a full rank lattice. Choose as a
fundamental cell for $\wG_1$ the set $\wK_1 := \delta_1\wedge\eta_1$ and let
$\wK_{11} := \{(x,y)\in\wK_1\st y \geq 1 - x\}$. Define as $\Omega_1$ the
set
\[
\Omega_1 := \wK_1\setminus\wK_{11} \cup (\wK_{11} + \delta_1 - 2\eta_1).
\]
Since $ \delta - 2\eta = (-\frac{1}{2},\frac{1}{2})^\top\in\wG_1$,
$\Omega_1$ is $\wG_1$-congruent to $\wK_1$.
\par
Now let $C_1 := \left\{(x,y)\in C\st y \leq x - \frac{1}{2}\right\}$, and let
$\rho_{\tilde{\alpha},1}\in\wW$ and $\rho_\alpha\in\wW$ denote the
reflection about the hyperplane $H_{\tilde{\alpha},1}$,
respectively, $H_\alpha$. As \[
\rho_\alpha\circ\rho_{\tilde{\alpha},1} (C_1) = \wK_1 + \delta_1 - 2\eta_1,
\]
the set $C$ is also $\wW$-congruent to $\Omega$.

For a second example, we choose new vectors $\delta_2$ and $\eta_2$ as follows: $ \delta_2 := \frac{1}{4}\tilde{\alpha} = (\frac{1}{4},\frac{1}{4})^\top$ and $\eta_2 := (\alpha+\beta) = (1,0)^\top$. Let
\[
\wG_2 := \Z\delta_2\oplus\Z\eta_2.
\]
Then, as above, $\wG_2\cap\wW=\wG_2\cap\Gamma \supseteq 4\wG_2$, and thus $\wG_2\cap\wW=\wG_2\cap\Gamma$ is a full rank lattice. As a fundamental cell for $\wG_2$ select the set $\wK_2 := \delta_2\wedge\eta_2$ and let
$\wK_{21} := \{(x,y)\in\wK_2\st y \geq 1 - x\}$. Define $\Omega_2$ to be the set
\[
\Omega_2 := \wK_2\setminus\wK_{21} \cup (\wK_{21} - 2\delta_2).
\]
As $2 \delta_2 = (\frac{1}{2},\frac{1}{2})^\top\in\wG_2$, $\Omega_2$ is $\wG_2$-congruent to $\wK_2$.
\par
Now set $C_2 := \left\{(x,y)\in C\st y \leq \frac{1}{4}\right\}$, and denote by $\rho_{\beta}\in\wW$ the reflection about the hyperplane $H_{\beta}$. Then,
\[
C = \wK_{2} \cup \rho_\beta (\wK_{21} - 2\delta_2) = C_2 \cup \rho_\beta (\wK_{21} - 2\delta_2),
\]
showing that $C$ is also $\wW$-congruent to $\Omega_2$.
\begin{figure}[h]
\begin{center}
\includegraphics[width=4.5cm,height=4.5cm]{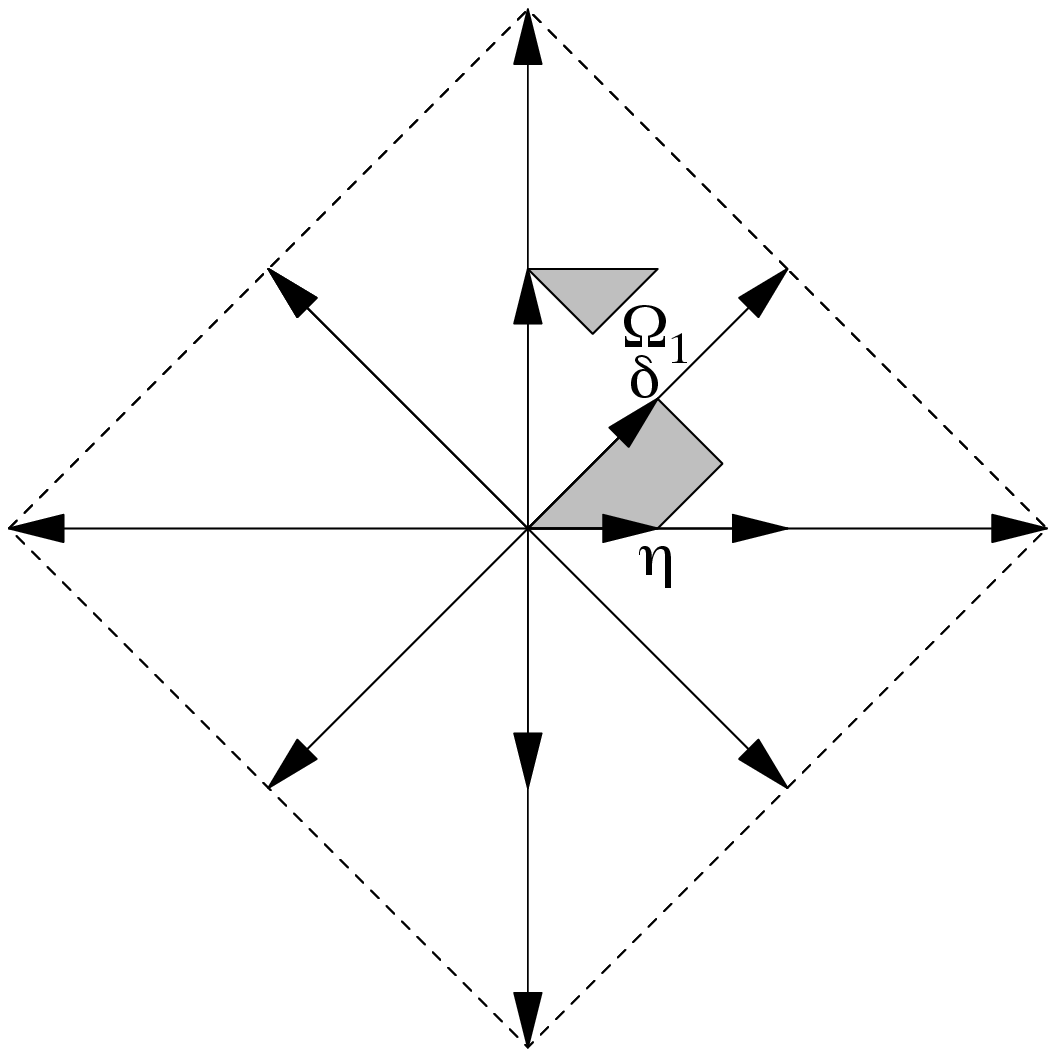} \hspace*{1cm}
\includegraphics[width=4.5cm,height=4.5cm]{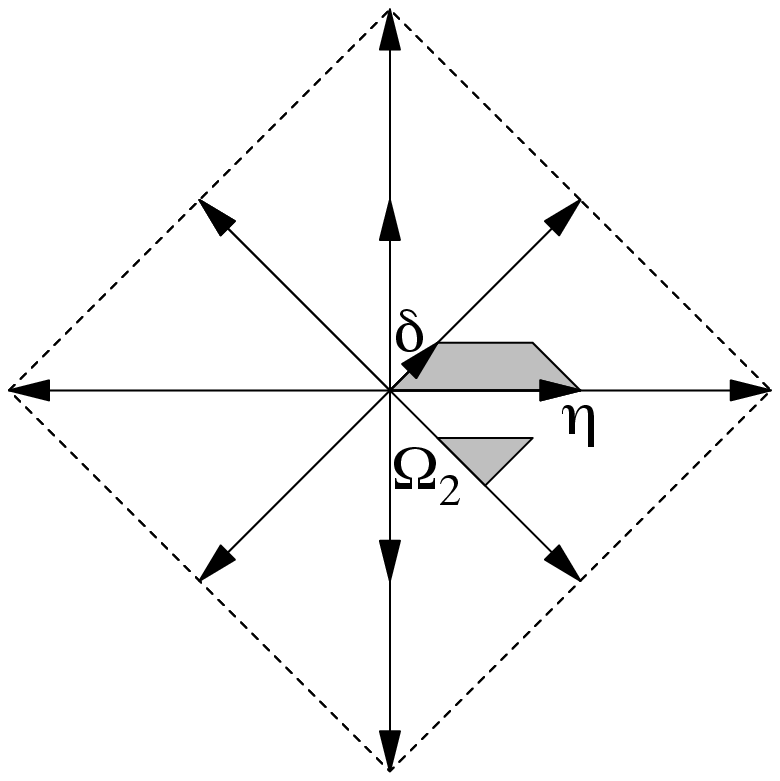} \caption{The sets $\Omega_1$ and $\Omega_2$ for
$B_2$.}\label{f1b}
\end{center}
\end{figure}
\subsection{Root System $G_2$}
\noindent In this final section we discuss the root system $G_2$
following the same line as in the last section.

As in Section \ref{s-1} let $V\subset\R^3$ be the hyperplane given
by the equation $x+y+z=0$ which we identify with $\R^2$ as before.

\begin{sumlist}
\item Root system $R$: Set of all $e\in L \cap V$ with length $\sqrt{2}$ or $\sqrt{6}$. These are: $\pm (e_1 - e_2)$,
$\pm (e_1 - e_3)$, $\pm (e_2 - e_3)$, $\pm (2 e_1 - e_2 - e_3)$,
$\pm (2 e_2 - e_1 - e_3)$, $\pm (2 e_13 - e_1 - e_2)$.

\item Positive roots: $\alpha$, $\beta$, $\alpha+\beta$, $2\alpha+\beta$, $3\alpha+\beta$, $3\alpha + 2\beta$.

\item Basis for $\R^2$: $\alpha := e_1 - e_2$ and $\beta := -2 e_1 + e_2 + e_3$.

\item Highest root: $\tilde{\alpha} = 3\alpha + 2\beta$.

\item Coroot system $R\,\check{}$:  $\pm (e_1 - e_2)$, $\pm (e_1 - e_3)$, $\pm (e_2 - e_3)$, $\pm \frac{1}{3}(2 e_1 - e_2 - e_3)$, $\pm \frac{1}{3}(2 e_2 - e_1 - e_3)$, $\pm \frac{1}{3}(2 e_3 - e_1 - e_2)$. Note that the short roots $\pm\alpha$, $\pm\alpha\pm\beta$, and $\pm 2\alpha\pm\beta$ are equal to their coroots, whereas the long roots have coroots 1/3 of their lengths.

\item As in case $A_2$ above, the angle between the roots $\alpha$ and $\beta$ can be computed to be $5\pi/6$ and a coordinate system $(v_1,v_2)$ in $V$ is introduced in the same fashion as above. In this coordinate system, $\alpha = (\sqrt{2},0)^\top$, $\beta = (-3\sqrt{2}/2,\sqrt{6}/2)^\top$, and $\tilde{\alpha} = 3\alpha + 2\beta = (0,\sqrt{6})^\top$.

\item In this case the Weyl group has order $|\cW|=12$.

\item Lattice generated by the coroots: $\Gamma =  \Z (\alpha + \frac{1}{3}\beta)\oplus  \Z (\alpha + \frac{2}{3}\beta)$.

\item Walls of the fundamental domain $C$ are given by the hyperplanes $H_\alpha: v_1 = 0$, $H_\beta: \sqrt{3} v_1 - v_2 = 0$, and $H_{\tilde{\alpha},1}: v_2 = 1/\sqrt{6}$.

\item Fundamental domain for $\widetilde{\mathcal{W}}$: $C$ is a triangle with vertices $(0,0)$, $(0,1/\sqrt{6})$, and $(\sqrt{2}/6,1/\sqrt{6})$, hence has angles $\pi/6$, $\pi/2$, $\pi/3$.

\item Area of $K$ and $C$ are: \[ |K| = \frac{\sqrt{6}}{3}\cdot\frac{\sqrt{2}}{2} =
\frac{\sqrt{3}}{3}\quad\textrm{and}\qquad |C| =
\frac{1}{2}\cdot\frac{\sqrt{2}}{6}\cdot\frac{1}{\sqrt{6}} =
\frac{\sqrt{3}}{36}. \]
\end{sumlist}
\par
\begin{figure}[h]
\begin{center}
\includegraphics[width=2in,height=2in]{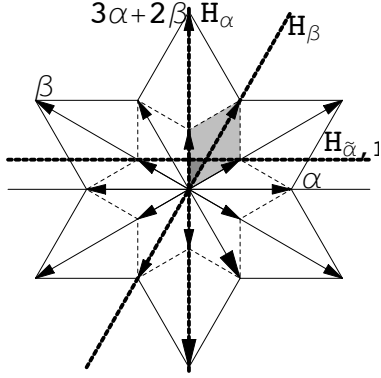} \caption{The geometry of $G_2$.}\label{f3}
\end{center}
\end{figure}
\par
Let $\delta_1 := \frac{1}{6}(3\alpha+2\beta) = (0,\sqrt{6}/6)^\top$, $\eta_1 := \frac{1}{6}(2\alpha+\beta) = (\sqrt{2}/12,\sqrt{6}/12)^\top$, and define $\wG_1 := \Z\delta_1\oplus\Z\eta_1$. Note that $6\wG_1\subseteq\Gamma$ and, thus, $\wG_1\cap\wW=\wG_1\cap\Gamma$ is a full rank lattice. Denote by $\wK_1$ the fundamental cell $\eta_1\wedge\delta_1$ of $\wG_1$, and set $\wK_{11} := \{(u,v)\in \wK\st v \geq \frac{\sqrt{6}}{6}\}$. Now define $\Omega_1$ as
\[
\Omega_1 := \wK_1\setminus\wK_{11}\cup (\wK_{11} +\delta_1 - 2\eta_1).
\]
Since $\delta_1 - 2\eta_1 = (-\frac{\sqrt{2}}{6},0)^\top\in\wG_1$, the set $\Omega_1$ is $\wG_1$-congruent to $\wK_1$. \par
Now let $C_1 := \left\{(u,v)\in C\st u \geq\frac{\sqrt{2}}{12}\right\}$, and let $\rho_{\tilde{\alpha},1}\in\wW$ and $\rho_\alpha\in\wW$ denote the reflection about the hyperplane $H_{\tilde{\alpha},1}$, respectively, $H_\alpha$. As
\[
\rho_\alpha\circ\rho_{\tilde{\alpha},1} (C_1) = \wK_1 + \delta_1 -2\eta_1,
\]
the set $C$ is also $\wW$-congruent to $\Omega_1$.

The second example is obtained by setting $\delta_2 := \frac{1}{12}(3\alpha+2\beta) = (0,\sqrt{6}/12)^\top$, $\eta_2 := \frac{1}{3}(2\alpha+\beta) = (\sqrt{2}/6,\sqrt{6}/6)^\top$, and defining $\wG_2 := \Z\delta_2\oplus\Z\eta_2$. Note that $6\wG_2\subseteq\Gamma$ and, thus, $\wG_2\cap\wW=\wG_2\cap\Gamma$ is a full rank lattice. Denote by $\wK_2$ the fundamental cell $\eta_2\wedge\delta_2$ of $\wG_2$, and let $\wK_{21} := \{(u,v)\in \wK\st v \geq 1\}$. Now define $\Omega_2$ as
\[
\Omega_2 := \wK_2\setminus\wK_{21}\cup (\wK_{21} - 2\delta_2 - \eta_2).
\]
Since $-2\delta_2 - \eta_2 \in\wG_2$, the set $\Omega_2$ is $\wG_2$-congruent to $\wK_2$.
\par
Let $C_2 := \left\{(u,v)\in C\st v \geq \sqrt{6} u + \frac{\sqrt{6}}{12}\right\}$, and let $\rho_\alpha\in\wW$ denote the reflection about the hyperplane $H_\alpha$. Then
\[
C = K_{21} \cup \rho_\alpha (\wK_{21} - 2\delta_2 - \eta_2) = C_2 \cup \rho_\alpha (\wK_{21} - 2\delta_2 - \eta_2),
\]
and the set $C$ is $\wW$-congruent to $\Omega_2$.
\begin{figure}[h] \begin{center}
\includegraphics[width=5cm,height=4cm]{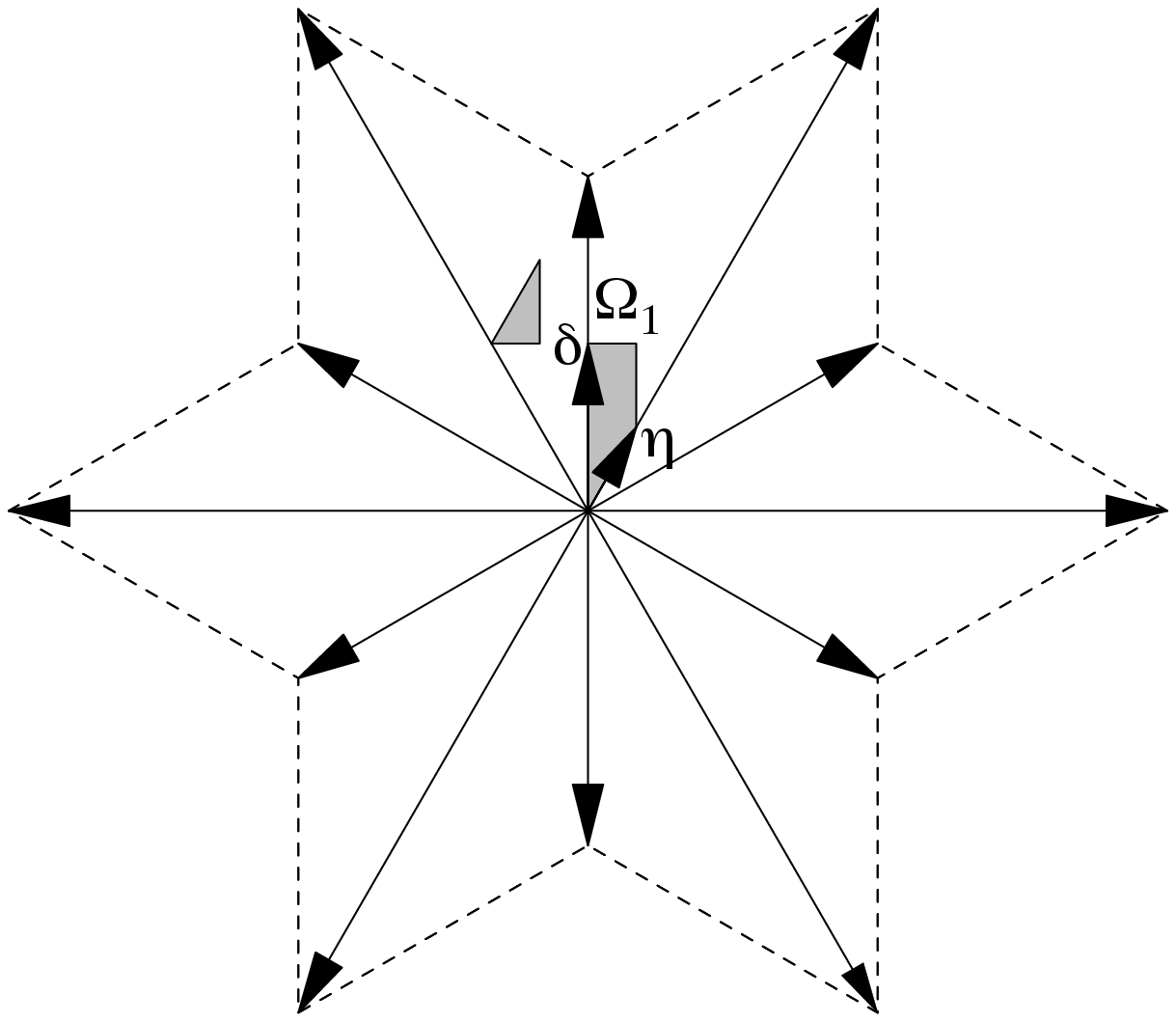} \hspace*{1cm}
\includegraphics[width=5cm,height=4cm]{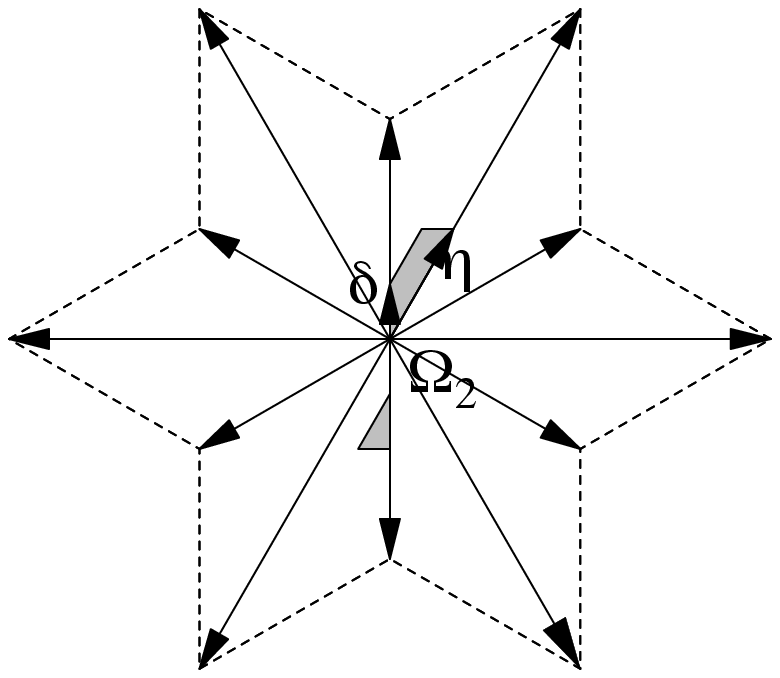}
\caption{The sets $\Omega_1$ and $\Omega_2$ for
$G_2$.}\label{f1c} \end{center} \end{figure}

The proof of Theorem 2.1 is complete.
\end{proof}

\begin{remark}
The set $\Omega_2$ for the root systems $B_2$ and $G_2$ was found by a student, F. Drechsler, while working on his diploma thesis \cite{Dr07} in mathematics at the Technische Universit\"{a}t M\"{u}nchen under the supervision of the second author.
\end{remark}

Combining Proposition 2.1  and Theorem 2.1, we have our main
theorem:

\begin{theorem} Let $A$ be an arbitrary $2\times2$ real expansive matrix, let $C$ be a foldable figure in $\R^2$ and let $\wW = \mathcal{W}\ltimes \Gamma$ be the associated affine Weyl group. Then,
there exists a full rank lattice $\wG$ and a measurable set $W$ that is a three-way tiling set for $\wG$,  $\wW$, and $\cD$. Hence, $W$ is simultaneously a $(\cD$,$\wG)$-dilation-translation wavelet
set and a $(\cD$,$\wW)$-dilation-reflection wavelet set.  So, in
particular, $(\wW, \wG, \cD)$ is an allowable triple.
\end{theorem}

\section{Problems}

Motivated by our results, we pose some new problems:
\medskip

\noindent
PROBLEM $2$.  If a triple $(\wW, \wG, \cD)$ has the property that
the foldable figure for $\wW$ and the lattice cell (parallelopiped)
for $\wG$ have the same volume, is the triple allowable?

An essential subproblem of this, which relates only translation and
reflection (and not dilation) is:
\medskip

\noindent
SUBPROBLEM $2a$.  If $\wW$ is an affine Weyl group of a foldable
figure in $\R^n$, and $\wG$ is the translation group of a full rank
lattice in $\R^n$, and if the measure of the foldable figure is the
same as the measure of a tiling cell for the lattice, does there
always exist a measurable set $K$ which tiles for both $\wW$ and
$\wG$? (This is a subproblem of Problem A because if $(\wW, \wG,
\cD)$ is allowable then $\wW$ and $\wG$ have a common tile by
definition.)

A natural strengthening of the hypothesis leads to the following
subproblem, which is intimately connected to the dual
translation-dilation/reflection-dilation wavelet theory of
\cite{lm06}.
\medskip

\noindent
SUBPROBLEM $2b$.  Let $\wW$ and $\wG$ be as in Subproblem $2a$,
and suppose in addition that the intersection of the two groups
contains the translation group of some full rank lattice. Does there
always exist a measurable set $K$ which tiles for both $\wW$ and
$\wG$? If such a dual reflection/translation tiling set exists, then
given any expansive matrix dilation $A$, an application of
Proposition 8.1 of \cite{lm06} yields a (perhaps different) set $W$
which is a three-way-tiler for $(\wW, \wG, \cD)$ for $\cD=\{A^n\mid
n\in \Z\}$. So a solution of this subproblem would in itself extend
our work significantly.
\bibliographystyle{amsalpha}
\end{document}